\documentclass[a4paper]{paper}

\usepackage{amsmath, amsfonts, amsthm}
\usepackage{xcolor}
\usepackage{graphicx}
\usepackage{caption}
\usepackage{datetime}
\usepackage{mathtools}

\newtheorem{thm}{Theorem}[section]
\newtheorem{lemma}[thm]{Lemma}
\theoremstyle{definition}
\newtheorem{remark}[thm]{Remark}
\newtheorem{assumption}[thm]{Assumption}

\usepackage{tikz}

\usetikzlibrary{calc}
\newcommand{\SpecimenInPlane}{\Pi_\text{in}}
\newcommand{\SpecimenOutPlane}{\Pi_\text{out}}

\newcommand{\LensPlane}{\Pi_\text{lens}}
\newcommand{\FocalPlane}{\Pi_\text{focal}}
\newcommand{\ImagePlane}{\Pi_\text{image}}

\newcommand{\incoming}{\text{in}}

\newcommand{\out}{\text{out}}

\newcommand{\wave}{u}

\newcommand{\TotalOp}{\mathcal{T}}

\newcommand{\ScatteringOp}{\TotalOp^\text{sc}}

\newcommand{\waveIn}{\wave^\text{in}}
\newcommand{\waveOut}{\wave^\text{out}}

\newcommand{\R}{\mathbb{R}}
\newcommand{\C}{\mathbb{C}}

\renewcommand{\Re}{\text{Re}}
\renewcommand{\Im}{\text{Im}}

\newcommand{\Intensity}{\mathcal{I}}
\newcommand{\RayOp}{\mathcal{P}}
\newcommand{\ScatteringPot}{f}
\newcommand{\ScatteringPotRe}{\ScatteringPot}
\newcommand{\ScatteringPotIm}{\ScatteringPot^\Im}
\newcommand{\ScatteringPotC}{\ScatteringPot^\C}

\newcommand{\Psf}{\text{PSF}}

\newcommand{\PropOp}{\mathcal{U}}

\newcommand{\Aperture}{A}
\newcommand{\FourierTransform}{\mathcal{F}}
\newcommand{\defocus}{{\Delta z}}

\mathtoolsset{showonlyrefs}

\title{Formal uniqueness in Ewald sphere corrected single particle analysis}
\subtitle{\today}
\author{P\"ar Kurlberg\thanks{P.K. acknowledges support by the Swedish Research
		Council (2016-03701, 2020-04036).} \and Gustav Zickert\thanks{G.Z. acknowledges support by the Swedish Research
		Council (2016-03701) and by the Swedish Foundation of Strategic Research under Grant AM13-0049. \newline Corresponding author. Email: gzickert@gmail.com}}

\begin{document}

\maketitle

\abstract{}
In single particle analysis (SPA), the task is to recover the scattering potential of a macromolecular structure from cryo-electron microscope images of many copies of the structure in unknown orientations. The idealized, noise-free SPA inverse problem has been shown to be uniquely solvable --- up to hand --- when the forward model is based on the ray transform. More accurate forward models take the non-zero curvature of the Ewald sphere into account.

We analyze an Ewald sphere corrected forward model for SPA and use the
diffraction slice theorem to prove that the corresponding inverse
problem is uniquely solvable, {\em including the hand of the
  structure}.

\section{Introduction}
A central problem in biology is to determine the structure of biological macromolecules. Single particle analysis (SPA) can be used to reveal the structure of molecules for which the traditional techniques of X-ray crystallography
and nuclear magnetic resonance spectroscopy have proved hard to use \cite{agard2014single, venien2017cryo, peplow2017cryo}.

SPA is an imaging technique in which a large number of ideally identical ``particles'' (e.g. protein molecules) in different relative orientations are imaged by a transmission electron microscope. The inverse problem of recovering the 3D structure of the particle from the resulting images is similar to other tomographic inverse problems such as X-ray computerized tomography, but what makes SPA rather special is that the projection directions are unknown. 

Let $\text{SE}(3)$ and $\text{SO}(3)$ denote the special Euclidean and special orthogonal groups in three dimensions, respectively. We use the convention that $\text{SE}(3)$ acts first by rotation and then by translation. More precisely,  $\text{SE}(3)$ and $\text{SO}(3)$ act on a function $f\colon \R^3 \to \R$ via  
\begin{align}
R\cdot f(x) &:= f(R^{-1}(x)), \quad R \in \text{SO}(3) \\
(R,c)\cdot f(x) &:= f(R^{-1}(x-c)), \quad (R,c) \in \text{SE}(3) 
\end{align}
If the center of mass of $f$ is at the origin, which we will assume in this text, then $c$ equals the center of mass of $(R,c)\cdot f$.

In the absence of noise, SPA data are often modelled as \cite{CGPW15, bendory2020single}: 
\begin{align}\label{eq:FullRayModel}
g_j = \Psf*\RayOp((R_j,c_j)\cdot f),\quad j\in J,
\end{align}
where $f\colon \R^3 \to \R_{\geq 0}$ represents the unknown particle, $J$ is some index set, $g_j\colon \R^{2} \to \R$ are the data, $(R_j, c_j) \in \text{SE}(3)$ are unknown rigid body motions, $\RayOp$ is a single-angle ray transform:
\begin{align}
\RayOp(f)(x_1, x_2) := \int_\R f(x_1,x_2, x_3)dx_3,
\end{align} 
and $\Psf$, the point spread function, is a given convolution kernel.

Clearly, SPA data from \eqref{eq:FullRayModel} can not distinguish $f$ from elements in the SE(3)-orbit of $f$. In fact we may only ever hope to recover the E(3)-orbit of $f$. In more practical terms, information of the \emph{hand} of $f$ (i.e. the orientation of $f$) is missing in data as in \eqref{eq:FullRayModel}. To see why this is the case, let $O := \text{diag}(1,1,-1)\in \text{O}(3)$, $R\in \text{SO}(3)$ and $R_O = ORO\in \text{SO}(3)$. Then  
\begin{align}
\RayOp(R_O\cdot(O\cdot f)) = \RayOp((R_OO)\cdot f) = \RayOp((OR)\cdot f) = \RayOp(O\cdot(R\cdot f)) = \RayOp(R\cdot f), 
\end{align}
where we used that $\RayOp(O\cdot f) = \RayOp(f)$. Hence, SPA data corresponding to $f$ and $\{R_j\}_{j\in J}$, is identical to data corresponding to $O\cdot f$ and $\{(R_j)_O\}_{j\in J}$.

Apart from these obstructions to uniqueness, the SPA inverse problem of recovering the structure $f$ from data given by \eqref{eq:FullRayModel} has been shown to be uniquely solvable using several different techniques. The method of common lines \cite{ lslss1986determination, van1987angular} is based on the fact that the Fourier transform of any two projection images will agree on at least one central line. For a generic triplet of projections this enables recovery of their mutual orientations (up to hand). In Kam's method \cite{kam1980reconstruction, bandeira2017estimation} one computes correlations of the Fourier transform of the structure using averaged correlations of the Fourier transforms of the projection images. Uniqueness was studied from a probabilistic point of view in \cite{panaretos2009random}, where the author proved that a pair of distinct objects induce mutually singular probability measures in the data space. Finally, there is the class of methods to which this paper belongs. It is based on relating moments of data to moments of the particle \cite{goncharov1988integral, salzman1990method, BB00, LY07, Lam08}. While it has not been widely used in practical reconstructions, it is useful for rigorously analyzing uniqueness.

The aforementioned studies all assume that data are originating from line integrals as in \eqref{eq:FullRayModel}. This \emph{projection assumption} corresponds to the assumption of a flat Ewald sphere, which limits the reconstruction quality at high resolutions \cite{glaeser2019good}. More accurate image formation models that better account for the wave properties of the imaging electron therefore take the curvature of the Ewald sphere into account \cite{derosier2000correction, FO08, russo2018ewald, leong2010correcting,wolf2006ewald, jensen2000defocus,wan2004full,kazantsev2010fully,voortman2014ctf}. In frequency space, these models provide data on half-spheres, rather than on planes\footnote{Some models \cite{wan2004full, voortman2014ctf} yield Fourier space data on paraboloids rather than on spheres.}.

Prior work on structure recovery based on Ewald sphere corrected models has mainly focused on the development of algorithms. However, to the best of our knowledge, there does not exist any explicit precise statement of uniqueness in the literature.

In this paper we focus on the noise-free and continuous setting. (The
latter means that we do not introduce any discretization of the
detector). In this setting we provide a mathematically rigorous proof
that accounting for a non-zero curvature of the Ewald sphere renders
the SPA inverse problem \emph{generically uniquely solvable}.

Our theorem treats a continuum limit where all SO(3) elements are
realized in the data (this assumption can be relaxed: the proof is
valid with minor modifications if the set of rotations in data is only
assumed to be a {\em countable dense subset of SO(3)}). The theorem also
assumes that all particles are located at the same distance from the
detector, that the real- and imaginary parts are related by a non-zero
constant amplitude contrast ratio, that $f$ has compact support and
that the first few moments of $f$ satisfy some generically valid
conditions (c.f. Assumption \ref{ass:Assym}). However, our result is
valid for any defocus and we assume neither that orientations are
known a priori, nor that the particle belongs to a certain symmetry
group, nor that multiple projections were recorded with varying
defocus. 

Since the intersection of two half-spheres is a circular segment, one could, perhaps, envision a uniqueness result based on finding the common curve associated to each pair of projection images. Orientation estimation based on such common curves has been considered in X-ray free electron laser imaging \cite{bortel2011common}, see also \cite{riedel2017geometric} for a discussion of the common curve problem in the context of SPA. 

The basic strategy here is to use the method of moments. Since the model we work with is more amenable for analysis in frequency space, we compute moments by differentiating the Fourier transform of data. A key element is to find the second order moments of $f$ by looking at extremal values of the second order moments of data, a technique which was also used in \cite{salzman1990method}.

The paper is organized as follows. In section \ref{sec:ForwardModel} we define our Ewald sphere corrected forward model. In section \ref{sec:Thm} we state and prove our uniqueness theorem. Finally, in section \ref{sec:Conclusion} we give a conclusion. 
\section{Forward model}\label{sec:ForwardModel}

We start out with the assumption that the unknown structure is described by a complex-valued and compactly supported scattering
 potential $\ScatteringPotC \in L^2(\R^3; \C) $ whose real- and
 imaginary part respectively determine the elastic- and inelastic
 scattering properties of the specimen. As mentioned in the
 introduction we assume that there exists a {\em known} constant $Q\in \R_{> 0}$, called the amplitude contrast ratio, such that real- and imaginary parts $\ScatteringPotRe$ and $\ScatteringPotIm$ of $\ScatteringPotC$ are related via
\begin{align}\label{eq:Paganini}
\ScatteringPotIm =Q\ScatteringPotRe.
\end{align}

The incoming electron is a spatial monochromatic plane wave $\waveIn(x) := e^{-ikx_3}$ with wave number $k>0$ travelling in the direction of $-e_3$. Its interaction  with the specimen is modelled by the Schr{\"o}dinger equation. Stationary solutions to the latter are given as solutions to a Helmholtz equation. In biological applications the specimen is typically weakly scattering and in this case the Born approximation, which linearizes the dependence on the scattering potential of the (approximate) solution to the Helmholtz equation, is applicable \cite{FO08, kohr2011fast}.
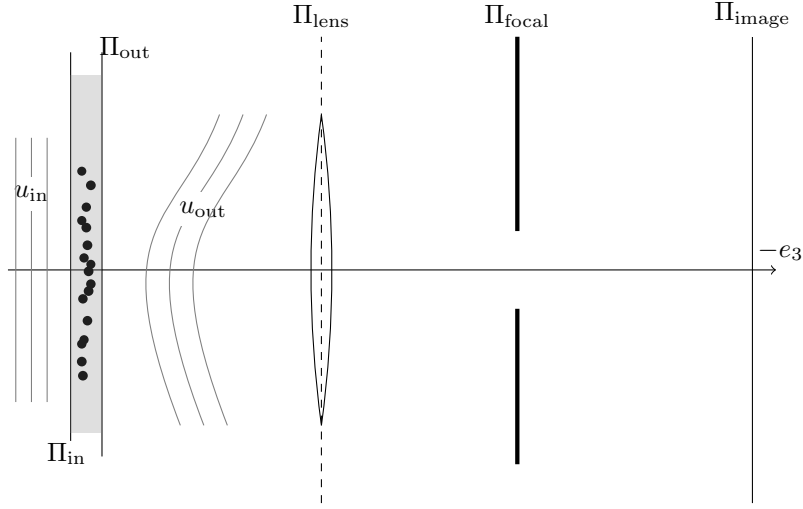
\begin{figure}[hbt!]
\centering
\newlength\imagewidth
\setlength{\imagewidth}{\textwidth}
\begin{tikzpicture}[scale=.85]
\coordinate (OpticalAxisStart) at (0,0);
\coordinate (OpticalAxisEnd) at (0.98\imagewidth,0); 
\coordinate (ObjectPlaneTop) at (0.1\imagewidth,0.3\imagewidth); 
\coordinate (ObjectPlaneBottom) at (0.1\imagewidth,-0.3\imagewidth); 
\coordinate (LensPlaneTop) at (0.4\imagewidth,0.3\imagewidth);
\coordinate (LensPlaneBottom) at (0.4\imagewidth,-0.3\imagewidth); 
\coordinate (ImagePlaneTop) at (0.95\imagewidth,0.3\imagewidth);
\coordinate (ImagePlaneBottom) at (0.95\imagewidth,-0.3\imagewidth);
\coordinate (LensNode) at (0.4\imagewidth,-0.2\imagewidth);
\coordinate (FocalPlaneTop) at (0.65\imagewidth,0.3\imagewidth);
\coordinate (FocalPlaneBottom) at (0.65\imagewidth,-0.25\imagewidth); 
\coordinate (ApertureTop) at (0.65\imagewidth,0.05\imagewidth);
\coordinate (ApertureBottom) at (0.65\imagewidth,-0.05\imagewidth); 
\coordinate (SpecimenCenter) at (0.1\imagewidth,0);
\coordinate (SpecimenInTop) at (0.08\imagewidth,0.28\imagewidth);
\coordinate (SpecimenInBottom) at (0.08\imagewidth,-0.22\imagewidth);
\coordinate (SpecimenInCenter) at (0.08\imagewidth,0);
\coordinate (SpecimenOutTop) at (0.12\imagewidth,0.28\imagewidth);
\coordinate (SpecimenOutBottom) at (0.12\imagewidth,-0.24\imagewidth);
\coordinate (SpecimenOutCenter) at (0.12\imagewidth,0);
\coordinate (SlabTop) at (0.12\imagewidth,0.25\imagewidth);
\coordinate (SlabBottom) at (0.08\imagewidth,-0.21\imagewidth);
\coordinate (ROITop) at (0.12\imagewidth,0.05\imagewidth);
\coordinate (ROIBottom) at (0.08\imagewidth,-0.05\imagewidth);
\coordinate (ElectronInTop) at (0.03\imagewidth,-0.17\imagewidth);
\coordinate (ElectronInBottom) at (0.03\imagewidth,0.17\imagewidth);
\coordinate (ScatteredElectronBottom) at (0.25\imagewidth,-0.2\imagewidth);
\coordinate (ScatteredElectronTop) at (0.3\imagewidth,0.2\imagewidth);
\draw[->] (OpticalAxisStart) -- (OpticalAxisEnd);
\node [above] at (OpticalAxisEnd) {{$\ -e_3$}}; 
\foreach \ypos in {-150, -130, -105, -99, -72, - 41, -30, -20, -2, 8, 17, 35, 60, 70, 89, 120, 140}
  {\pgfmathrandominteger{\xpos}{-4}{4}
    \pgfmathrandominteger{\radius}{70}{75} 
    \fill[rotate around={0:(SpecimenCenter)}]
             (0.1\imagewidth+0.5*\xpos,0.011*\ypos) circle (0.001*\radius);
   };
\newcommand*{\TiltAngleInPlot}{0}%
\draw[rotate around={\TiltAngleInPlot:(SpecimenInCenter)}] 
  ([rotate around={\TiltAngleInPlot:(SpecimenInCenter)}] SpecimenInBottom) -- 
  ([rotate around={\TiltAngleInPlot:(SpecimenInCenter)}] SpecimenInTop);
\node [below] at ($ ([rotate around={0:(SpecimenInCenter)}] SpecimenInBottom) + (-0.005\imagewidth,0.01\imagewidth) $) {{$\SpecimenInPlane$}}; 
\draw[rotate around={-10:(SpecimenOutCenter)}]
  ([rotate around={\TiltAngleInPlot:(SpecimenOutCenter)}] SpecimenOutBottom) -- 
  ([rotate around={\TiltAngleInPlot:(SpecimenOutCenter)}] SpecimenOutTop);
\node [above] at ($ ([rotate around={0:(SpecimenOutCenter)}] SpecimenOutTop) + (0.03\imagewidth,-0.02\imagewidth) $) {{$\SpecimenOutPlane$}}; 
\draw [draw=none,fill=gray,fill opacity=.25,rotate around={\TiltAngleInPlot:(SpecimenCenter)}]%
  ([rotate around={\TiltAngleInPlot:(SpecimenInCenter)}] SlabBottom) rectangle
  ([rotate around={\TiltAngleInPlot:(SpecimenInCenter)}] SlabTop);   
\node [above] at ($ ([rotate around={-10:(SpecimenCenter)}] SpecimenOutTop) - (0.03\imagewidth,0.07\imagewidth) $) {{}};
\draw[dashed] (LensPlaneBottom) -- (LensPlaneTop);
\node [above, fill=white] at (LensPlaneTop) {{$\LensPlane$}}; 
\draw (LensNode) arc (-7.65:7.65:1.5\imagewidth);
\begin{scope}[yscale=1,xscale=-1, xshift=-0.8\imagewidth]
  \draw (LensNode) arc (-7.65:7.65:1.5\imagewidth);
\end{scope}
\draw [ultra thick] (FocalPlaneBottom) -- (ApertureBottom);
\draw [ultra thick] (ApertureTop) -- (FocalPlaneTop);
\node [above, fill=white] at (FocalPlaneTop) {{$\FocalPlane$}}; 
\draw (ImagePlaneBottom) -- (ImagePlaneTop);
\node [above, fill=white] at (ImagePlaneTop) {{$\ImagePlane$}}; 
\draw[gray] (ElectronInBottom) -- (ElectronInTop);
\node [below, fill=white] at ($ (ElectronInBottom) - (0,0.05\imagewidth) $) {{$u_{\incoming}$}}; 
\draw[gray] ($ (ElectronInBottom) - (0.02\imagewidth,0) $) -- ($ (ElectronInTop) - (0.02\imagewidth,0) $);
\draw[gray] ($ (ElectronInBottom) + (0.02\imagewidth,0) $) -- ($ (ElectronInTop) + (0.02\imagewidth,0) $);
\draw[gray] (ScatteredElectronBottom) .. controls (1.8,0.8) and (3,0.7) .. (ScatteredElectronTop);
\node [below, fill=white] at ($ (ScatteredElectronTop) - (0.051\imagewidth,0.1\imagewidth) $) {{$u_{\out}$}}; 
\pgftransformxshift{-0.03\imagewidth}
\draw[gray] ($ (ScatteredElectronBottom) - (0.03\imagewidth,0) $) .. controls (1.8,0.8) and (3,0.7) .. ($ (ScatteredElectronTop) - (0.03\imagewidth,0) $);
\pgftransformxshift{0.06\imagewidth}
\draw[gray] ($ (ScatteredElectronBottom) + (0.03\imagewidth,0) $) .. controls (1.8,0.8) and (3,0.7) .. ($ (ScatteredElectronTop) + (0.03\imagewidth,0) $);
\end{tikzpicture}
\caption{\footnotesize Schematic representation of the SPA set-up. For the uniqueness theorem we assume that the relative shifts along the optical axis are small enough to be considered equal to zero.}
\label{fig:OpticalSetUp}
\end{figure}
Under the Born approximation, which we use in this paper, the interaction of $\waveIn$ with $\ScatteringPotC$ results in an outgoing wave $u^\text{out}$ at the specimen exit plane $\SpecimenOutPlane := e_3^\perp \cong \R^2$.

The scattering operator $\ScatteringOp$ maps $\ScatteringPotC$ to the corresponding $\waveOut$:
\begin{align}
\ScatteringOp\left(\ScatteringPotC\right)(x) := u^\text{out}(x) = 1+k^{-1}\PropOp\left(\ScatteringPotC\right)(x) ,\quad  x\in \R^2
\end{align}
where
$\PropOp$ is the so-called propagation operator \cite{NW01, FO08} that is defined as
\begin{align}
\PropOp\left(\ScatteringPotC\right) := \frac{k}{\waveIn}G_k *\left(\waveIn\ScatteringPotC\right)\bigg\rvert_{e_3^\perp}.
\end{align}
Here $*$ denotes convolution on $\R^3$ and 
\begin{align}
G_k(x) := \frac{1}{4\pi}\frac{e^{ik|x|}}{|x|}, \quad x\in \R^3 \setminus \{0\}.
\end{align}

The outgoing wave then passes through the optics of the electron microscope. The effect of the optics is here modelled by the convolution\footnote{Formally, we define convolution by point-wise multiplication in frequency space.}
\begin{align}
\waveOut \mapsto h * \waveOut, 
\end{align}
where the point spread function $h$ is defined by its Fourier transform:
\begin{align}
& \FourierTransform\left[h\right](\xi) := A(\xi)e^{i\chi(\xi)}, \quad \chi(\xi) := a|\xi|^2 + b|\xi|^4.
\end{align}
In the above, the aperture function $A$ is the indicator function of some ball and
\begin{align}
a &:= \frac{\defocus}{2k} \in \R \\
b &:= -\frac{C_s}{4k^3} \in \R,  
\end{align}
where $C_s$ is a constant that encodes spherical aberrations of the optics and $\defocus$ is the defocus. We have assumed perfectly coherent illumination for simplicity.
Finally, the wave forms an intensity distribution in the detector plane. This leads to the following total non-linear model for the measured intensity:
\begin{align}
	\Intensity_\text{non-lin}\left(\ScatteringPotC\right)(x) := \left\lvert h * \left\{ 1+k^{-1}\PropOp\left(\ScatteringPotC\right)\right\}(x)\right\rvert^2, x\in \R^2.
\end{align}

Next,  we assume that the quadratic term
is negligible and study instead the linearized model
	\begin{align}\label{eq:TotalForwardLinear}
& \Intensity_0(\ScatteringPotC) :=  k^{-1} h *\PropOp\left(\ScatteringPotC\right) +   k^{-1}\overline{h}*\overline{\PropOp\left(\ScatteringPotC\right)}.
\end{align}
With assumption \eqref{eq:Paganini} in place we define our final intensity operator $\Intensity$ that acts directly on a real-valued potential via \begin{align}
\Intensity(\ScatteringPotRe) := \Intensity_0((1+iQ)\ScatteringPotRe) = k^{-1}(1+iQ) h *\PropOp\left(\ScatteringPotRe\right) +   k^{-1}\overline{(1+iQ)}\overline{h}*\overline{\PropOp\left(\ScatteringPotRe\right)}.
\end{align}

\begin{remark}
It can be shown \cite{FO08, NW01} that, point-wise in frequency space, one has
\begin{align}
&2 \cdot \PropOp\left(\ScatteringPotRe\right) \xrightarrow[]{k\to \infty} i\RayOp\left(\ScatteringPotRe\right), 
\end{align}
so for large  $k$ the Born approximation reproduces the classical model: 
\begin{align}  
&	\Intensity\left(\ScatteringPotRe\right) \approx k^{-1}\Re\left((i-Q) h\right) * \RayOp\left(\ScatteringPotRe\right) = -k^{-1}\Psf_Q * \RayOp\left(\ScatteringPotRe\right), \\
& \FourierTransform\left[\Psf_Q\right](\xi) :=  1_{B(0,r)}(\xi)\left[Q\cos\chi(\xi) + \sin\chi(\xi)\right].
\end{align}
\end{remark}
So far we have defined the forward model corresponding to the electron microscope image of a single particle. In single particle analysis, the complete data is a collection
\begin{align}
\left\lbrace\Intensity\left((R_j,c_j)\cdot\ScatteringPotRe\right)\ \middle \vert\ j\in J \right\rbrace
\end{align}
 of rotated and translated copies of $\ScatteringPotRe$, where the $R_j$'s and $c_j$'s are unknown and $J$ is some index set. A schematic drawing of the set-up is given in Fig. \ref{fig:OpticalSetUp}.
\section{A uniqueness result}\label{sec:Thm}
\subsection{Some preparatory definitions, lemmas and assumptions}
The moments of a function $f\colon \R^n\to \C$ are defined as 
\begin{align}
m_\alpha f &:= \int_{\R^n}x^\alpha fdx,\quad \alpha \in \mathbb{N}^n,
\end{align}
where $\alpha = (\alpha_1,\dots,\alpha_n) \in \mathbb{N}^n$ is a multi-index, and $x^\alpha := \prod_ix_i^{\alpha_i}$. The order of a multi-index $\alpha$ is the natural number $|\alpha|$ defined by $|\alpha| := \sum_i \alpha_i$. We define the matrix 
\begin{align}
\Lambda(\ScatteringPotRe) := \left(\begin{matrix}
m_{200}\ScatteringPotRe & m_{110}\ScatteringPotRe & m_{101}\ScatteringPotRe \\
m_{110}\ScatteringPotRe & m_{020}\ScatteringPotRe & m_{011}\ScatteringPotRe \\
m_{101}\ScatteringPotRe & m_{011}\ScatteringPotRe & m_{002}\ScatteringPotRe 
\end{matrix}\right).
\end{align}

For $f\colon \R^n\to \C$ we use the following version of the Fourier transform:
\begin{align}
\FourierTransform(f)(\xi):=\hat f(\xi) &:= \int_{\R^n}f(x)e^{-i\xi\cdot x}dx.
\end{align}
In order to carry out our analysis in frequency space, we need the following two basic lemmas.
\begin{lemma}\label{lem:momDeriv}
	\begin{align}
	(-i)^{|\alpha|}m_\alpha f = \partial^\alpha \hat f (0),\quad \alpha\in \mathbb{N}^n.
	\end{align}
\end{lemma}
\begin{lemma}\label{lem:SE3Fourier}
	\begin{align}
	\FourierTransform\left[(R,c)\cdot f)\right] (\xi) = e^{-ic\cdot\xi}\left(R\cdot \hat f\right)(\xi),\quad (R,c) \in \text{SE}(3).
	\end{align}
\end{lemma}
The space of functions that may serve as real-valued scattering potentials is
\begin{align}
X_0 := \left\lbrace \ScatteringPotRe\in L^2(\R^3)\ \middle \vert\ \ScatteringPotRe \text{ has compact support and } \ScatteringPotRe\geq 0 \right \rbrace.  
\end{align}
We will require the following mild generic asymmetry assumptions on the scattering potential.
\begin{assumption}\label{ass:Assym}
\begin{enumerate} 
	\item[]
	\item $\Lambda(\ScatteringPotRe)$ has three distinct (real) eigenvalues.
	\item In a coordinate-system where $\Lambda(\ScatteringPotRe)$
          is diagonal, neither of the third order moments $m_{300}$ and $m_{210}$ vanish.
\end{enumerate}
\end{assumption}
\noindent Finally, let 
\begin{align}
X := \left\lbrace \ScatteringPotRe\in X_0 \ \middle \vert\ \ScatteringPotRe \text{ satisifies Assumption } \ref{ass:Assym} \right\rbrace.
\end{align}
\subsection{Theorem statement and proof}
\begin{thm}
Consider the data
\begin{align}
y := \left\lbrace\Intensity\left((R_j,c_j)\cdot\ScatteringPotRe\right)\ \middle \vert\ (R_j, c_j) \in \text{SE}(3), j\in J \right\rbrace,
\end{align}
for some index set $J$. If $\{R_j\ \vert\ j\in J\} = \text{SO}(3)$
 and there exits a known positive constant $c_0$ such that $\forall j\colon c_j \cdot e_3 = c_0$, then $y$ determines $\ScatteringPotRe\in X$ uniquely up to a rigid body motion.  
\end{thm}
\begin{remark}
  As already mentioned, a minor modification of the argument shows
  that the conclusion holds for any countable set $J$ such that
  $\{R_j\ \vert\ j\in J\}$ is dense in SO(3). We leave the details to
  the interested reader.
  \end{remark}
\begin{proof}
To fix the position and pose of $\ScatteringPotRe$, we assume that
\begin{enumerate}
	\item$f$ has center of mass in the origin, which implies that $\nabla \widehat{\ScatteringPotRe} (0) =0$. 
	\item $\Lambda(\ScatteringPotRe)=\text{diag}(m_{200}, m_{020}, m_{002})$, with $m_{002} < m_{020} < m_{200}$. 
	 \item $m_{300}, m_{210} > 0$.
\end{enumerate}
Since $f$ is compactly supported, its Fourier transform is analytic and therefore uniquely determined by its Taylor series coefficients. The basic idea of the proof is to show that (up to a rigid body motion) these coefficients can be uniquely recovered from $y$ via a system of polynomial equations whose unknowns are given by the coefficients along with the particle translations and rotations.
The uniqueness of the solution to the system of equations is established by induction over the order of the Taylor coefficients. Once a coefficient order is fixed, we are able to reduce our system to the one that corresponds to the ray transform, and the latter is known to have a unique solution.
Unless otherwise stated, in the proof below, ``\dots is known'' should be interpreted as ``\dots is uniquely determined by $y$ up to a rigid body motion''.

\subsubsection*{Step 1: Fourier transform data and apply the diffraction slice theorem}
For $\xi$ small enough so that the aperture function satisfies $\Aperture(\xi)=1$, the Fourier transform of \eqref{eq:TotalForwardLinear} is given by
\begin{align}
\FourierTransform\left[k\Intensity_0(\ScatteringPotC)\right] =  e^{i\chi}\FourierTransform\left[\PropOp\left(\ScatteringPotC\right)\right]   + e^{-i\chi}\FourierTransform\left[\overline{\PropOp\left(\ScatteringPotC\right)}\right] \label{eq:ForwardModelFourierExp}.
\end{align}
Assuming a constant amplitude contrast ratio (i.e. $\ScatteringPotC = (1+iQ)\ScatteringPotRe$), the RHS of \eqref{eq:ForwardModelFourierExp} can be expressed as
\begin{align}
\omega_1 \FourierTransform\left[\PropOp\left(\ScatteringPotRe\right)\right] + \overline{\omega_1} \FourierTransform\left[\overline{\PropOp\left(\ScatteringPotRe\right)}\right],
\end{align}
where $\omega_1 := (1+iQ)e^{i\chi}$. We further rewrite this using the identity $\overline{\FourierTransform\left[h\right]}(\xi) = \FourierTransform\left[\bar h\right](-\xi)$ as 
\begin{align}\label{eq:ForwardModelFourierExpPaganini}
\omega_1 \FourierTransform\left[\PropOp\left(\ScatteringPotRe\right)\right] + \overline{\omega_1} \overline{\FourierTransform\left[\PropOp\left(\ScatteringPotRe\right)\right]}(-\cdot).
\end{align}
The diffraction slice theorem \cite{NW01} provides us with the following expression for $\FourierTransform\left[\PropOp\left(\ScatteringPotRe\right)\right]$:
\begin{align}
\FourierTransform\left[\PropOp\left(\ScatteringPotRe\right)\right](\xi) = \phi(\xi)\widehat{\ScatteringPotRe}   \left(\gamma^+ (\xi)\right),
\end{align}
where \begin{align}
\phi(\xi) := \frac{i}{2} \frac{k}{k-\gamma_3(\xi)}
\end{align} and $\gamma^+$ is the map that lifts the plane to a half-sphere:
\begin{align}
\gamma^+ (\xi) &:= \left(\xi_1, \xi_2, \gamma_3(\xi)\right) \\
\gamma_3(\xi) &:= k - \sqrt{k^2 - |\xi|^2}
\end{align}
Next, again using the identity $\overline{\FourierTransform\left[h\right]}(\xi) = \FourierTransform\left[\bar h\right](-\xi)$, we compute: 
\begin{align}
\overline{\FourierTransform\left[\PropOp\left(\ScatteringPotRe\right)\right]}(-\xi) &=  \bar\phi(-\xi) \overline{\widehat{\ScatteringPotRe} \left(\gamma^+ (-\xi)\right)} 
=\bar\phi(\xi) \widehat{\ScatteringPotRe}   \left(-\gamma^+ (-\xi)\right) 
=\bar\phi(\xi) \widehat{\ScatteringPotRe}   \left(\gamma^-(\xi)\right),
\end{align}
where $\gamma^-$ maps the plane to a half-sphere that extends downwards. More precisely, it is defined as 
\begin{align}
\gamma^-(\xi) := \left(\xi_1, \xi_2, -\gamma_3(\xi)\right).
\end{align}
Hence, from
\eqref{eq:ForwardModelFourierExpPaganini} it follows that
\begin{align}
\FourierTransform\left[\Intensity_0\left(\ScatteringPotC\right)\right] &=\omega_2 \widehat{\ScatteringPotRe} \circ \gamma^+ + \overline {\omega_2} \widehat{\ScatteringPotRe} \circ \gamma^-,
\end{align}
where 
\begin{align}
\omega_2 := \frac{1}{k}\omega_1\phi = \frac{1}{2}(Q-i)e^{i\chi} \frac{1}{\gamma_3-k}.
\end{align}
Hence we obtain the Fourier space model
\begin{align}
2\left(\gamma_3-k\right)\FourierTransform\left[\Intensity_0\left(\ScatteringPotC\right)\right] &=ze^{i\chi} \widehat{\ScatteringPotRe} \circ \gamma^+ +\bar{z}e^{-i\chi} \widehat{\ScatteringPotRe} \circ \gamma^-,
\end{align}
where 
\begin{align}
z := (Q-i).
\end{align}
\subsubsection*{Step 2: Handle in-plane translations}
For $(R, c) \in \text{SE}(3)$ we introduce the notation
\begin{align}
h^{(1)}_{R,c}(\xi) &:= 2\left(\gamma_3-k\right)\FourierTransform\left[\Intensity((R,c)\cdot\ScatteringPotRe)\right](\xi) \\
&= ze^{i\chi}\FourierTransform\left[(R,c)\cdot \ScatteringPotRe\right] \circ \gamma^+ + \bar z e^{-i\chi}\FourierTransform\left[(R,c)\cdot \ScatteringPotRe\right]\circ \gamma^- \\
&= ze^{i\chi}e^{-ic\cdot\gamma^+}\left(R\cdot\widehat{\ScatteringPotRe} \right)\circ \gamma^+ + \bar z e^{-i\chi}e^{-ic\cdot\gamma^-} \left(R\cdot\widehat{\ScatteringPotRe} \right)\circ \gamma^- \\
&= e^{-i(c_1, c_2)\cdot \xi}\left[ze^{i\chi}e^{-ic_0\gamma_3}\left(R\cdot\widehat{\ScatteringPotRe} \right)\circ \gamma^+ + \bar z e^{-i\chi}e^{+ic_0\gamma_3} \left(R\cdot\widehat{\ScatteringPotRe} \right)\circ \gamma^-\right].  \label{eq:c_tilde_factored_out}
\end{align}
(To simplify the notation we suppress the
dependence on $\xi$ and simply write
$\gamma^{\pm}$, $\gamma_{3}$, $\xi$ and $\widehat{\ScatteringPotRe}$.) 

Set $a_\alpha := \partial^\alpha \widehat{\ScatteringPotRe} (0)/\alpha!$. Since $\ScatteringPotRe$ has compact support, it follows from the Paley-Wiener theorem (see e.g. \cite[Theorem 7.1.14]{Hor03}) that $\widehat{\ScatteringPotRe}$ can be extended to an entire function on $\C^n$, and in particular the following Taylor expansions in $\xi \in \R^2$ are well-defined. 
\begin{align}
\gamma_3 &= \frac{1}{2k}|\xi|^2 + \mathcal{O}\left(|\xi|^4\right) \\
e^{i\left(\pm\chi(\xi) \mp c_0\gamma_3\right) } &= 1 \pm i\left(a  -\frac{c_0}{2k}\right)|\xi|^2 + \mathcal{O}\left(|\xi|^4\right) \\
\widehat{\ScatteringPotRe} &= \widehat{\ScatteringPotRe}(0) + \sum_{|\alpha| = 2,3}a_\alpha\xi^\alpha + \mathcal{O}\left(|\xi|^4\right)
\end{align}
Hence \eqref{eq:c_tilde_factored_out} implies that
\begin{align}\label{eq:FirstOrderTaylor}
h^{(1)}_{R,c}(\xi) = (z + \overline{z})\widehat{\ScatteringPotRe}(0) - i(z + \overline{z})\widehat{\ScatteringPotRe}(0)(c_1, c_2)\cdot \xi + \mathcal{O}\left(|\xi|^2\right).
\end{align}
We note that since $Q$ is assumed to be known, $\widehat{\ScatteringPotRe}(0)$ is known from the zeroth order term of \eqref{eq:FirstOrderTaylor}. Therefore,  $(c_1, c_2)$ can be read of from the first order term of \eqref{eq:FirstOrderTaylor}. Thus $\left((c_j)_1,(c_j)_2\right)_{j\in J}$ is known, hence the set 
\begin{align}
H_1 := \left\{h^{(2)}_{R}(\xi)\mid R\in \text{SO}(3)\right\}
\end{align}
is known, where
\begin{align}
h^{(2)}_{R}(\xi) &:=ze^{i\chi}e^{-ic_0\gamma_3}\left(R\cdot\widehat{\ScatteringPotRe} \right)\circ \gamma^+ + \bar z e^{-i\chi}e^{+ic_0\gamma_3} \left(R\cdot\widehat{\ScatteringPotRe} \right)\circ \gamma^-. 
\end{align} 
\subsubsection*{Step 3: Extract a one-parameter subset of data}
Note that
\begin{align}
h^{(2)}_{R}(\xi)  &= z\left(1 + i\left(a  -\frac{c_0}{2k}\right)|\xi|^2\right)\left(\widehat{\ScatteringPotRe}(0) + \sum_{|\alpha| = 2,3}a_\alpha(R^{-1}\gamma^+(\xi))^\alpha \right) \\
&+ \bar z\left(1 - i\left(a  -\frac{c_0}{2k}\right)|\xi|^2\right)\left(\widehat{\ScatteringPotRe}(0) + \sum_{|\alpha| = 2,3}a_\alpha(R^{-1}\gamma^-(\xi))^\alpha \right) + \\
&+ \mathcal{O}\left(|\xi|^4\right) \\
&= (z+\bar z)\widehat{\ScatteringPotRe} (0) + C|\xi|^2 \\
&+ z\left(a_{200}(R^{-1}\gamma^+(\xi))_1^2 + a_{020}(R^{-1}\gamma^+(\xi))_2^2 + a_{002}(R^{-1}\gamma^+(\xi))_3^2\right) \\
&+ \bar z\left(a_{200}(R^{-1}\gamma^-(\xi))_1^2 + a_{020}(R^{-1}\gamma^-(\xi))_2^2 + a_{002}(R^{-1}\gamma^-(\xi))_3^2\right) \\
&+ z\sum_{|\alpha| =3}a_\alpha(R^{-1}\gamma^+(\xi))^\alpha + \bar z\sum_{|\alpha| =3}a_\alpha(R^{-1}\gamma^-(\xi))^\alpha + \mathcal{O}\left(|\xi|^4\right),
\end{align} 
where $C:= 2 \widehat{\ScatteringPotRe}(0)\text{Re}\left(zi\left(a-\frac{c_0}{2k}\right)\right)$

Since $\widehat{\ScatteringPotRe}(0)$ is known, $C$ is known. Now if
\begin{align}
R^{-1} = \left(\begin{matrix}
R_{11} & R_{12} & R_{13} \\
R_{21} & R_{22} & R_{23} \\
R_{31} & R_{32} & R_{33}
\end{matrix}\right),
\end{align}
then
\begin{align}
\left[R^{-1}\gamma^\pm(\xi)\right]_j^2 &=\left(R_{j1}\xi_1 + R_{j2}\xi_2 \pm R_{j3}\frac{1}{2k}(\xi_1^2 + \xi_2^2)\right)^2
+ \mathcal{O}(|\xi|^3) \\
&=R_{j1}^2\xi_1^2 + 2R_{j1}R_{j2}\xi_1\xi_2 + R_{j2}^2\xi_2^2 
+ \mathcal{O}(|\xi|^3),
\end{align}
and therefore the $\xi_1^2$-term in the expansion of $h^{(2)}_{R}$ is
\begin{align}
(z+\bar z)\left(a_{200}R_{11}^2+a_{020}R_{21}^2+a_{002}R_{31}^2\right) + C.
\end{align}
Note that Lemma \ref{lem:momDeriv} and $m_{200}f > m_{020}f > m_{002}f$ together imply that $a_{200} < a_{020} < a_{002}$, so the $\xi_1^2$-term has a minimal value of $(z+\bar z) a_{200}+C$, which is achieved exactly when $R_{11} = \pm 1, R_{21}=R_{31}=0$. The maximal value is $(z+\bar z)a_{002}+C$. Minimizing the $\xi_1^2$-term thus forces $R^{-1}$ to have the form
\begin{align}
R^{-1} = R_{S_1, \theta}^{-1} := \left(\begin{matrix}
S_1 & 0 & 0 \\
0 & \cos \theta & -\sin \theta \\
0 & \sin \theta & \cos \theta
\end{matrix}\right),
\end{align}
for some $\theta\in [0, 2\pi)$ and $S_1 \in \{-1, 1\}$.

This means that the set 
\begin{align}
H_2 := \left\{h^{(2)}_{R_{S_1},\theta }(\xi)\mid \theta\in [0, 2\pi), S_1 \in \{-1, 1\}\right\}
\end{align}
is known.

Now we take a closer look at the third-order expansion of $ h^{(2)}_{R_{S_1, \theta}}$:
\begin{align}
h^{(2)}_{R_{S_1, \theta}}(\xi) &= (z+\bar z)\widehat{\ScatteringPotRe}(0) + \left((z+\bar z) a_{200}+C\right)\xi_1^2 \\
&+ \left((z+\bar z) a_{020}  \cos^2(\theta) +(z+\bar z) a_{002}  \sin^2(\theta) +C\right)\xi_2^2 \\
&+ \biggl((z+\bar z)a_{210}   \cos(\theta) +(z+\bar z)  a_{201}  \sin(\theta) \\
&\qquad + \frac{z - \bar z}{k}(a_{002} - a_{020}) \cos(\theta)\sin(\theta)\biggr) \xi_1^2\xi_2 \\
&+(z+\bar z) S_1a_{300}  \xi_1^3 + \beta_1(S_1,\theta)\xi_1\xi_2^2 + \beta_2(S_1,\theta)\xi_2^3 + \mathcal{O}\left(|\xi|^4\right),
\end{align}
for some functions $\beta_1, \beta_2 \colon \{-1, 1\} \times [0 , 2\pi) \to i\R$. Observe that the sign of the imaginary part of the $\xi_1^3$-coefficient equals $S_1$. Hence, with $R_\theta := R_{1, \theta}$, the following set is known:
\begin{align}
H_3 :=\left\{h^{(2)}_{R_\theta }(\xi)\mid \theta\in [0, 2\pi)\right\}.
\end{align}

\subsubsection*{Step 4: Compute the second order moments of $f$}
In what follows, unless otherwise is stated, ``coefficient'' refers to
a Taylor coefficient of $h^{(2)}_{R_{\theta}}$.  Note first that
$a_{200}$ is determined by the $\xi_1^2$ coefficient. Likewise,
$a_{020}$ and $a_{002}$ are known, since the minimal and maximal
values (recall that the $\xi_1^2$ coefficient is minimal) of the
$\xi_2^2$ coefficients equals $a_{020}+ C$ and $a_{002}+C$,
respectively (note that $C$ was determined in step 3).

\subsubsection*{Step 5: Compute some third order moments of $f$}
The $\xi_2^2$-coefficient is minimized exactly when
$\cos\left(\theta\right) = \pm 1$ and $\sin\left(\theta\right) =
0$. Hence, an examination of the $\xi_1^2\xi_2$-coefficient shows that
$|a_{210}|$ is known. Thus $a_{210}$ is known, since the sign of
$a_{210}$ was fixed in the beginning of this proof. Similarly, by 
maximizing the $\xi_2^2$-coefficient, we can determine $|a_{201}|$.

\subsubsection*{Step 6: Recover the hand of $\ScatteringPotRe$ along with some rotations}
The $\xi_2^2$-coefficient determines $\cos^2\left(\theta\right)$,
hence each $\theta$ is determined up to a most four
possibilities. More precisely, each $\theta$ is determined up to sign
and a shift by $\pi$. 

We will now see that all small $\theta$ may be determined. There exists a known positive constant $\epsilon
= \epsilon(z, a_{020}, a_{002}, a_{210}, |a_{201}|)  < \pi/2$ such
that if $\min\{|\theta|, |\theta-\pi|\} < \epsilon$, then the sign of the imaginary part of
the $\xi_1^2\xi_2$-coefficient\footnote{Note that the
  $\xi_1^2\xi_2$-coefficient is purely imaginary.} equals
the sign of its pure cosine term.
Let 
\begin{align}
\sigma(\theta) &:= (z+\bar z) a_{020}  \cos^2(\theta) +(z+\bar z) a_{002}  \sin^2(\theta) +C
\end{align}
and for $0 < \theta < \epsilon$ consider the known set
\begin{align}
y_1(\theta) := \Bigl\{h\in H_3  \ \Big\vert\ &  \text{The }\xi_2^2\text{-coefficent of } h \text{ equals }\sigma(\theta) \text{ and the}\\
& \text{imag. part of the }\xi_1^2\xi_2\text{-coefficent of } h \text{ is negative }\Bigr\}.
\end{align}
Recall that $m_{210} > 0$ by assumption 3 in the beginning of this proof. It then follows from Lemma \ref{lem:momDeriv} and the discussion above that for any $\varphi$ and any $\theta \in (0,\epsilon)$:
\begin{align}
h^{(2)}_{R_{\varphi}} \in y_1(\theta) \Rightarrow \cos(\varphi) > 0, 
\end{align}
  so we have that
\begin{align}
y_1(\theta) = \left\{h^{(2)}_{R_{\theta}}, h^{(2)}_{R_{-\theta}}\right\}, 0 < \theta < \epsilon.
\end{align}
To remove the final sign-ambiguity, we study the
$\xi_1^2\xi_2$-coefficient. Note that
$|\sin(\theta)| = \sin(|\theta|)$ is known for $|\theta| < \epsilon$,
since $|\theta|$ is known. We fix some $\theta_0 \in (0,\epsilon)$ and
further fix some 
$h_0\in y_1(\theta_0)$. Let $h_{21}$ denote the
$\xi_1^2\xi_2$-coefficient of $h_0$ and let
\begin{align}
D := \frac{\left\lvert h_{21} - \cos(\theta_0)(z+\bar z)a_{210}\right\rvert}{|\sin(\theta_0)|}.
\end{align}
Note that
\begin{align}
D &= |A + B|, \\
A &= (z+\bar z)a_{201},\\
B &= \frac{z - \bar z}{k}(a_{002} - a_{020}) \cos(\theta_0).
\end{align}
In the above $|A|, B$ and $D$ are known and $|A|, B\neq 0$\footnote{Recall that $Q>0$ by assumption, so $z+ \bar z \neq 0$.}. Thus $A$ is known, so $a_{201}$ is known.

For $\theta\in (0,\epsilon)$ and $h\in y_1(\theta)$, let $h_{21}(\theta)$ denote the $\xi_1^2\xi_2$-coefficient of $h$. Then
\begin{align}\label{eq:extractSine}
h_{21}(\theta) - \cos(\theta)(z+\bar z)a_{210} = \Bigl(\pm\sin(\theta)\Bigr) \left((z+\bar z)  a_{201} + \frac{z - \bar z}{k}(a_{002} - a_{020}) \cos(\theta)\right).
\end{align}
Hence $\sin(\theta)$ is known for $0<\theta < \epsilon$ (the right factor in the RHS of \eqref{eq:extractSine} might vanish for some $\theta$: in this case we shrink $\epsilon$ accordingly.). Since $\cos(\theta)$ was already known, $\theta $ is known. More precisely, for $\theta \in (0, \epsilon)$ we now know
\begin{align}
h^{(2)}_{R_\theta}(\xi) &=ze^{i\chi}e^{-ic_0\gamma_3}\left(R_\theta\cdot\widehat{\ScatteringPotRe} \right)\circ \gamma^+ + \bar z e^{-i\chi}e^{+ic_0\gamma_3} \left(R_\theta\cdot\widehat{\ScatteringPotRe} \right)\circ \gamma^-.
\end{align}
\subsubsection*{Step 7: Recover all moments of $f$} 
Since 
\begin{align}
\left[R_\theta^{-1}\gamma^\pm(\xi)\right]_1^i\left[R_\theta^{-1}\gamma^\pm(\xi)\right]_2^j\left[R_\theta^{-1}\gamma^\pm(\xi)\right]_3^k 
\end{align}
equals
\begin{align}
&\xi_1^i\left(\cos(\theta)\xi_2 + \mathcal{O}
\left(|\xi|^{2}\right)\right)^j\left(\sin(\theta)\xi_2+ \mathcal{O}
\left(|\xi|^{2}\right)\right)^k\\
&= \xi_1^i\cos^j(\theta)\sin^k(\theta)\xi_2^{j+k} + \mathcal{O}
\left(|\xi|^{i+j+k+1}\right),
\end{align}
it follows that the Taylor series for $h^{(2)}_{R_\theta}$ admits the decomposition
\begin{align}
(z+\bar z)^{-1}h^{(2)}_{R_\theta} &=  \sum_{i,j,k}a_{ijk} \xi_1^i\cos^j(\theta)\sin^k(\theta)\xi_2^{j+k}  + \sum_{r,s}b_{rs}(\theta) \xi_1^r\xi_2^s,
\end{align}
where $a_{ijk}$ denotes the $ijk$:th Taylor coefficient of $\widehat{\ScatteringPotRe}$ and $b_{rs}$ only depends on $a_{ijk}$ for $i,j,k$ such that $i+ j+k < r+s$.

Hence, if $c_{ij}$ is defined as the $\xi_1^i\xi_2^j$-coefficient in the known Taylor series of $(z+\bar z)^{-1}h^{(2)}_{R_\theta}$, then
\begin{align}\label{eq:Radon2DmomentEqs}
\sum_{k=0}^ja_{ik(j-k)}\cos^k(\theta)\sin^{j-k}(\theta) = c_{ij} - b_{ij}(\theta),
\end{align}
where $b_{ij}(\theta)$ depends only on $a_{rst}$ for $r,s,t$ such that $r+s+t < i+j$.

Fix a moment-order $m\in \mathbb{N}$ and assume that $a_{ijk}$ has been determined for all $i,j,k$ such that $i+j+k < m$. Moreover fix $i$ and $j$ such that $i+j=m$ and let $\theta_1,\dots,\theta_{j+1} \in (0,\epsilon)$ be distinct angles. We consider the linear system of $j+1$ equations in the $j+1$ variables $\left(a_{ik(j-k)}\right)_{k=0}^j$ that results from inserting those angles in \eqref{eq:Radon2DmomentEqs}:
\begin{align}\label{eq:Radon2DmomentEqFinite}
\sum_{k=0}^ja_{ik(j-k)}\cos^k(\theta_\ell)\sin^{j-k}(\theta_\ell) = c_{ij} - b_{ij}(\theta_\ell), \quad \ell=1,\dots, j+1.
\end{align}
The system of equations \eqref{eq:Radon2DmomentEqFinite} is known (c.f. \cite{milanfar1996moment}) to be uniquely solvable for the variables $\left(a_{ik(j-k)}\right)_{k=0}^j$.\footnote{Note that the right hand side of \eqref{eq:Radon2DmomentEqFinite} is known by assumption.}

Consequently, $a_{ijk}$ is known for all $i,j,k$ such that $i+j+k=m$, and by induction $a_{ijk}$ is known for any $i,j,k\in \mathbb{N}$. Since $\hat{\ScatteringPotRe}$ can be extended to an entire function, it is uniquely determined by its Taylor series coefficients $(a_{ijk})_{i,j,k=0}^\infty$. Hence $\hat \ScatteringPotRe$ is known, which implies that  $\ScatteringPotRe$ is known. This concludes the proof.
\end{proof}

\section{Conclusion}\label{sec:Conclusion}
Existing uniqueness results for SPA applies to a forward model that is based on the ray transform. A more accurate diffraction tomographic model that accounts for the curvature of the Ewald sphere is offered in the Born approximation framework. Although models based on the ray transform have been used successfully for SPA structure recovery, the attainable resolution is ultimately limited by the incorrect assumption of a flat Ewald sphere.

We extend the classical method of moments, which is based on relating moments of data to moments of the unknown 3D structure, from the ray transform model to a model based on the Born approximation. Through the diffraction slice theorem this allows us to prove that the SPA inverse problem based on an Ewald sphere corrected model is uniquely solvable, including the hand of the structure.

 We believe that several sharpenings and extensions of the theorem are
 possible, e.g. it seems likely that some of the assumptions (e.g. on
 the set of needed rotations and translations) can be relaxed.

\section*{Acknowledgments}
The authors thank Ozan \"Oktem for suggesting the idea to study uniqueness in SPA with Born approximation and for indispensable guidance and feedback. 

\bibliographystyle{alpha}
\bibliography{ref}
\end{document}